\documentclass[12pt]{article}

\usepackage{latexsym}
\usepackage{amsfonts}
\usepackage{amssymb}
\usepackage{amsmath}
\usepackage{mathrsfs}
\usepackage{bbm}
\usepackage{dsfont}
\usepackage{hyperref}

\newcommand{\be}{\begin{equation}}
\newcommand{\ee}{\end{equation}}
\newcommand{\bea}{\begin{eqnarray}}
\newcommand{\eea}{\end{eqnarray}}
\newcommand{\bean}{\begin{eqnarray*}}
\newcommand{\eean}{\end{eqnarray*}}
\newcommand{\brray}{\begin{array}}
\newcommand{\erray}{\end{array}}

\newtheorem{dfn}{Definition}[section]
\newtheorem{thm}[dfn]{Theorem}
\newtheorem{lmma}[dfn]{Lemma}
\newtheorem{ppsn}[dfn]{Proposition}
\newtheorem{crlre}[dfn]{Corollary}
\newtheorem{xmpl}[dfn]{Example}
\newtheorem{rmrk}[dfn]{Remark}

\newcommand{\bdfn}{\begin{dfn}\rm}
\newcommand{\bthm}{\begin{thm}}
\newcommand{\blmma}{\begin{lmma}}
\newcommand{\bppsn}{\begin{ppsn}}
\newcommand{\bcrlre}{\begin{crlre}}
\newcommand{\bxmpl}{\begin{xmpl}}
\newcommand{\brmrk}{\begin{rmrk}\rm}

\newcommand{\edfn}{\end{dfn}}
\newcommand{\ethm}{\end{thm}}
\newcommand{\elmma}{\end{lmma}}
\newcommand{\eppsn}{\end{ppsn}}
\newcommand{\ecrlre}{\end{crlre}}
\newcommand{\exmpl}{\end{xmpl}}
\newcommand{\ermrk}{\end{rmrk}}

\newcommand{\cla}{\mathcal{A}}

\newcommand{\clh}{\mathcal{H}}

\newcommand{\clg}{\mathcal{G}}

\author{ S. Sundar}
\title {Product systems associated to compound Poisson processes}

\begin{document}
\maketitle
\begin{abstract}
In this paper, we consider a simple test case of  multiparameter product systems that arise out of random measures. We associate a product system to a stationary Poisson process and a stationary compound Poisson process. We show that the resulting $E_0$-semigroups are CCR flows.
%We also exhibit a two parameter decomposable product system that does not admit any unit. 
 %This paper is an application of the results obtained in \cite{Injectivity}.
  \end{abstract}
\noindent {\bf AMS Classification No. :} {Primary 46L55; Secondary 60G55.}  

\section{Introduction}
%The theory of $E_0$-semigroups initiated by R.T. Powers and developed extensively by Arveson is now a well established field of research and is  more than three decades old. 
%For a comprehensive account of $1$-parameter $E_0$-semigroups, we refer the reader to \cite{Arveson}.  Recently the author in collaboration with others (\cite{Anbu_Sundar}, \cite{Murugan_Sundar_continuous})
%has studied $E_0$-semigroups in the multiparameter context. 

Let $P$ be a closed convex cone in a Euclidean space $\mathbb{R}^{d}$. We assume that $P$ spans $\mathbb{R}^{d}$ and $P$ contains no line.  An $E_0$-semigroup over $P$ is
a semigroup of normal unital $*$-endomorphisms, indexed by $P$, of the algebra of bounded operators of an infinite dimensional separable Hilbert space satisfying a natural continuity assumption. It was first established by Arveson (\cite{Arv_Fock4}, \cite{Arv}), in the one parameter context, that $E_0$-semigroups are in
one-one correspondence with product systems. This was recently extended to the multiparameter case in \cite{Murugan_Sundar_continuous}. Roughly speaking, a product system over $P$ is a measurable
field of separable Hilbert spaces carrying an associative multiplication rule. %It was proved in \cite{Murugan_Sundar_continuous} and \cite{Murugan_Sundar} that every $E_0$-semigroup over $P$ has an associated
%product system and every product system over $P$ is isomorphic to a product system associated to a unique $E_0$-semigroup (up to cocycle conjugacy). 

%Although the bijection between abstract product system and $E_0$-semigroups was established in \cite{Murugan_Sundar_continuous}, no example of an abstract product system was  given. Let us elucidate what we 
%mean a bit more. The product systems that we know come from $E_0$-semigroups (ex, CAR flows, CCR flows). That is one defines the $E_0$-semigroup and obtain the product system as a consequence. It is 
%desirable to have an example of an abstract product system so that the bijection established in \cite{Murugan_Sundar} have some meaning. The aim of this paper is to produce  such examples and
%identify the resulting $E_0$-semigroups. 

In the one parameter setting, examples of  product systems, which are in fact exotic, were constructed by Tsirelson in (\cite{Tsirelson}, \cite{Tsi}, \cite{Tsirelson_gauge}) and by Liebscher in \cite{Liebscher} by making use of  probabilistic models. The remarkable works 
of Tsirelson and that of Liebscher have amply demonstrated the rich interaction that exists between probability theory and the theory of $E_0$-semigroups. We expect similar things to take place in the multiparameter setting
too.  We explore a simple test case here. 

The simplest probabilistic models that give rise to product systems in the multiparameter setting are stationary compound Poisson processes. In this paper, we  define a product system corresponding to a 
stationary compound Poisson process. This is akin to associating product systems to  L\'{e}vy processes (\cite{Mark}) in the one parameter setting. We show that in the case of a stationary Poisson process
and in the case of a stationary compound Poisson process, the resulting $E_0$-semigroups  are CCR flows. 

In the Poisson case, this is expected due to the Wiener-It\^{o} Chaos decomposition
(See Chapter 18, \cite{Last}). However it is quite easy to exhibit an explicit isomorphism in the Poisson case and it does not require sophisticated knowledge of Poisson processes. We then use the results of the Poisson case and the computation of the local projective cocycles carried out in 
\cite{Injectivity} to identify the $E_0$-semigroup associated to a stationary compound Poisson process. 

 To keep the prerequistes on point processes to a minimum, in Section 2, we have collected the preliminaries on  Poisson processes. For completeness, we have included proofs of some results. Our exposition is based on \cite{Kingman} and \cite{Last}.
 % In Section 3, we  associate a product system to a stationary Poisson process and verify that it is indeed a product system. In Section 4, we make a small digression into decomposable product systems. We explain that a decomposable product system, in the multiparameter setting, admitting a unit is a CCR flow. This is the conceptual reason behind the fact that the product systems associated to  Poisson processes are CCR flows.  We illustrate with a two paramter example that unlike in the $1$-parameter case, it is not true in general that a decomposable product system admits a unit. In Section 5, we identify the product system associated to a stationary Poisson process and a stationary compound Poisson process. 
 We make the following assumptions throughout this paper. 
\begin{enumerate}
\item[(1)] The probability triples $(\Omega,\mathcal{F},\mathbb{P})$ that we consider are assumed to be complete.
\item[(2)] Let $(\Omega,\mathcal{F},\mathbb{P})$ be a probability space. All sub $\sigma$-algebras of $\mathcal{F}$ are assumed to be complete. Let $\Sigma \subset \mathcal{F}$ be a sub $\sigma$-algebra. Recall that $\Sigma$ is said to be complete if it satisfies the following. If $A \subset \Omega$ equals a set in $\Sigma$ up to measure zero then $A \in \Sigma$. 
\item[(3)] For random variables $X_i$, by the smallest $\sigma$-algebra generated by $X_i$, we mean the smallest complete $\sigma$-algebra which makes $X_i$'s measurable. 
\end{enumerate}

The author would like to thank Prof. Partha Sarathi Chakraborty for his suggestion to investigate the relation between point processes and $E_0$-semigroups. 
\section{Poisson processes}
Let $(\mathbb{X},\mathcal{B})$ be a measurable space. We assume that $\mathcal{B}$ is countably generated. Let $\overline{\mathbb{N}}=\mathbb{N} \cup \{\infty\}$. 
Denote the set of $\overline{\mathbb{N}}$-valued $\sigma$-finite measures by $\mathcal{N}(\mathbb{X})$. We endow $\mathcal{N}(\mathbb{X})$ with the smallest $\sigma$-algebra which makes the map \[\mathcal{N}(\mathbb{X}) \ni \mu \to \mu(B) \in \overline{\mathbb{N}}\] measurable for every $B \in \mathcal{B}$. 

A point process on $\mathbb{X}$ is a random element of $\mathcal{N}(\mathbb{X})$, i.e. a measurable mapping $\eta: \Omega \to \mathcal{N}(\mathbb{X})$ where $(\Omega,\mathcal{F},\mathbb{P})$ is an underlying probability space. Let $\eta$ be a point process on $\mathbb{X}$. For $\omega \in \Omega$, $\eta(\omega)$ is a measure on $\mathbb{X}$ and for $B \in \mathcal{B}$, the map $\omega \to \eta(\omega)(B)$ is a random variable. We denote the random variable $\omega \to \eta(\omega)(B)$ by $\eta(B)$. 

  Let $\lambda$ be a $\sigma$-finite measure on $\mathbb{X}$. A point process $\eta$ on $\mathbb{X}$ is called a Poisson process with intensity measure $\lambda$ if the following two conditions are satisfied.
  \begin{enumerate}
  \item[(1)] For $B \in \mathcal{B}$, $\eta(B)$ is a Poisson random variable with parameter $\lambda(B)$, i.e. for $k \in \mathbb{N}$, $\mathbb{P}(\eta(B)=k)=\frac{(\lambda(B))^{k}}{k!}e^{-\lambda(B)}$.  
  \item[(2)] The random variables $\eta(B_1),\eta(B_2),\cdots,\eta(B_n)$ are independent whenever  $(B_i)_{i=1}^{n}$ is a disjoint family of measurable subsets of $\mathbb{X}$. 
  \end{enumerate}
  For a $\sigma$-finite measure $\lambda$, there exists a unique (up to equality in distribution) Poisson process with intensity measure $\lambda$. 
    For a measure $\nu$ on $\mathbb{X}$ and an integrable  complex valued function $u$ on $\mathbb{X}$, we denote $\int u(x)\nu(dx)$ by $\nu(u)$. 

\begin{lmma}
\label{Master equation}
Let $\eta$ be a Poisson process on $\mathbb{X}$ with intensity measure $\lambda$. Let $u$ be a simple function on $\mathbb{X}$ which is integrable. Then 
\[
\mathbb{E}(e^{\eta(u)})=exp\Big(\int (e^{u(x)}-1)\lambda(dx)\Big).\]
\end{lmma}
\textit{Proof.} Write $\displaystyle u=\sum_{i=1}^{n}a_i1_{B_i}$ with $B_i$'s disjoint and $a_i \in \mathbb{C}\backslash\{0\}$. Note that $\lambda(B_i)<\infty$ for every $i$. Then $\displaystyle e^{\eta(u)}=\prod_{i=1}^{n}e^{a_i\eta(B_i)}$. Making use of the independence of $\eta(B_i)$ and the fact that $\eta(B_i)$ is a Poisson random variable with parameter $\lambda(B_i)$, calculate as follows to observe that 
\begin{align*}
\mathbb{E}(e^{\eta(u)})&=\prod_{i=1}^{n}\mathbb{E}(e^{a_i\eta(B_i)}) \\
&=\prod_{i=1}^{n}exp(\lambda(B_i)(e^{a_i}-1)) \\
&=exp(\sum_{i=1}^{n}(e^{a_i}-1)\lambda(B_i))\\
&=exp\Big(\int (e^{u(x)}-1)\lambda(dx)\Big).
\end{align*}
The proof is now complete. \hfill $\Box$

Let $(\Omega,\mathcal{F},\mathbb{P})$ be an underlying probability space realising the Poisson process $\eta$ with intensity measure $\lambda$. We can and will assume that $\mathcal{F}$ is the smallest $\sigma$-algebra which makes the family of random variables $\eta(B)$, $B \in \mathcal{B}$, measurable.

\begin{ppsn}
\label{total}
The linear span of $\{e^{\eta(u)}: \textrm{$u$ is a non-negative simple $L^{2}$-function}\}$, is dense in $L^{2}(\Omega,\mathcal{F},\mathbb{P})$. 
\end{ppsn}
See Lemma 18.4 of \cite{Last} for a proof. 
%\textit{Proof.} Let $\mathbb{G}$ be the real linear span of $\{e^{\eta(u)}: \textrm{$u$ is a real valued simple $L^{2}$-function}\}$. Note that $\mathbb{G}$ is closed under multiplication. Let $W$ be the set of all bounded real valued measurable functions on $\Omega$ which can approximated in $L^{2}$ by elements of $G$. Note that 
%\begin{enumerate}
%\item[(a)] $W$ is a vector space containing constant functions,
%\item[(b)] if $f_n$ is an increasing sequence of non-negative functions in $W$ such that $\sup_n\{f_n\}<\infty$ and $f=\lim_{n \to \infty} f_n$ then $f \in W$, and
%\item[(c)] $W$ is closed under uniform convergence.
%\end{enumerate}
%Hence $W$ contains all bounded $\sigma(\mathbb{G})$ measurable functions. Let $B$ be a measurable subset of $\mathbb{X}$ such that $\lambda(B)<\infty$. Note that the pointwise limit of $\frac{e^{t\eta(B)}-1}{t}$ as $t \to 0$ is $\eta(B)$. Hence $\eta(B)$ is $\sigma(\mathbb{G})$ measurable for every $B$ with finite measure. For an arbitrary measurable set $B$, write $B=\coprod_{n=1}^{\infty}B_n$ where $B_n$ has finite measure. Then $\eta(B)=\sum_{n=1}^{\infty}\eta(B_n)$. As a consequence $\eta(B)$ is measurable for every $B \in \mathcal{B}$. Hence $\mathcal{F} \subset \sigma(\mathbb{G})$ which implies that $\mathcal{F}=\sigma(\mathbb{G})$. Consequently, $W$ contains every real valued bounded measurable function on $\Omega$. From here it is routine to complete the proof. \hfill $\Box$

\begin{ppsn}
\label{separability}
The space $L^{2}(\Omega,\mathcal{F},\mathbb{P})$ is separable. 
\end{ppsn}
\textit{Proof.} Let $\cla \subset \mathcal{B}$ be an algebra such that $\cla$ is countable and $\cla$ generates $\mathcal{B}$. Enumerate the elements of $\cla$ as $A_1,A_2,\cdots$. Choose $B_n \in \mathcal{B}$ such that $\lambda(B_n)<\infty$ and $\coprod_{n=1}^{\infty}B_n=\mathbb{X}$. It suffices to prove that $\mathcal{F}$ is the smallest $\sigma$-algebra which makes the random variables $\eta(A_m \cap B_n)$ measurable. To that end, let $\mathcal{G}$ be the smallest $\sigma$-algebra which makes the random variables $\eta(A_m \cap B_n)$, $m,n \in \mathbb{N}$, measurable.  

For $n \in \mathbb{N}$, let $\clg_n:=\{B \in \mathcal{B}: B \subset B_n, \eta(B) \textrm{ is $\mathcal{G}$ measurable}\}$. It is immediate that $\clg_n$ is a monotone class and contains the algebra $\cla \cap B_n$. Hence $\clg_n$ coincides with the $\sigma$-algebra of subsets of $B_n$ generated by $\cla \cap B_n$. Consequently $\clg_n=\mathcal{B} \cap B_n$. Now for $B \in \mathcal{B}$, $\eta(B)=\sum_{n=1}^{\infty}\eta(B \cap B_n)$. Hence $\eta(B)$ is $\clg$-measurable for every $B \in \mathcal{B}$. Hence $\clg=\mathcal{F}$. This completes the proof. \hfill $\Box$

\section{Product systems associated to  stationary Poisson processes} 

In this section, we associate a product system to a stationary Poisson process. Let $P$ be a closed convex cone in $\mathbb{R}^{d}$ which we assume is spanning and pointed, i.e. $P-P=\mathbb{R}^{d}$ and $P \cap -P=\{0\}$. Denote the interior of $P$ by $Int(P)$. For $x,y \in \mathbb{R}^{d}$, we write $x \leq y$ if $y-x \in P$ and we write $x<y$ if $y-x \in Int(P)$. The cone $P$ is fixed for the reminder of this paper. The setting that we consider is as follows.

Let $(Y,\mathcal{B})$ be a measurable space on which $\mathbb{R}^{d}$ acts in a measurable fashion.  We use additive notation for the action. Assume that $\mathcal{B}$ is countably generated.   Let $\lambda$ be a $\sigma$-finite measure on $Y$ which is   $\mathbb{R}^{d}$ invariant. Consider a Poisson process $\eta$ on $Y$ with intensity measure $\lambda$. Then $\eta$ is stationary, i.e. for measurable subsets $B_1,B_2,\cdots,B_n$ of $Y$ and $x \in \mathbb{R}^{d}$,  $(\eta(B_1),\cdots,\eta(B_n))=(\eta(B_1+x),\eta(B_2+x),\cdots,\eta(B_n+x))$ in distribution.

Suppose that $X\in \mathcal{B}$ and is $P$-invariant, i.e. $x+a \in X$ whenever $x \in X$ and $a \in P$. We assume that the action of $P$ on $X$ is pure, i.e. $\displaystyle \bigcap_{a \in P}(X+a)=\emptyset$. Note that for $a,b \in P$ with $a\leq b$, $X+b \subset X+a$. Let $(\Omega,\mathcal{F},\mathbb{P})$ be an underlying probability space realising the Poisson procees $\eta$. For $a,b \in P$ with $a \leq b$, let $\mathcal{F}_{a,b}$ be the $\sigma$-algebra generated by $\{\eta(B): B \subset (X+a )\backslash( X+b), B \in \mathcal{B}\}$. For $a \in P$, set $\mathcal{F}_{a}:=\mathcal{F}_{0,a}$. 

Note that for $a \leq b \leq c$, $(X+a)\backslash(X+c)=((X+a)\backslash(X+b))\coprod ((X+b)\backslash (X+c))$. This together with the complete independence of the Poisson process implies that for $a \leq b \leq c$, $\mathcal{F}_{a,b}$ and $\mathcal{F}_{b,c}$ are independent and $\mathcal{F}_{a,c}$ is generated by $\{\mathcal{F}_{a,b}, \mathcal{F}_{b,c}\}$. 

The stationarity of the Poisson process implies that for every $c \in P$, there exists a unitary $S_{c}:L^{2}(\Omega,\mathcal{F}_{a,b},\mathbb{P}) \to L^{2}(\Omega, \mathcal{F}_{a+c,b+c}, \mathbb{P})$ such that 
\[
S_c(f(\eta(B_1),\eta(B_2),\cdots,\eta(B_n)))=f(\eta(B_1+c),\eta(B_2+c),\cdots,\eta(B_n+c)).\]

We are now in a position to define the product system given the above data. Let 
\begin{equation}
\label{product system}
E:=\{(a,f): a \in P, f \in L^{2}(\Omega, \mathcal{F}_{a},\mathbb{P})\}.\end{equation}
The first projection of $E$ onto $P$ is denoted by $p$. Note that by Lemma \ref{separability}, the fibres $E(a)$, for $a \in P$, are separable. The product structure on $E$ is defined as follows
\[
(a,f)(b,g)=(a+b,fS_{a}(g)).\]

It is routine to verify from discussions above that the algebraic requirements for $E$ to be a product system are met. We proceed towards the proof of the fact that $E$ is  a product system in the sense of Definition 9.2 of \cite{Sundar_Notes}.
A careful look at the axioms of Definition 9.2 of \cite{Sundar_Notes} reveals that it suffices to prove the following. 
\begin{enumerate}
\item[(1)] The set $E$ is a Borel subset of $P \times L^{2}(\Omega, \mathcal{F}, \mathbb{P})$. 
\item[(2)] The multiplication $E \times E \ni (a,f) \times (b,g) \to (a+b, fS_a(g))$ is measurable.
\item[(3)] The pair $(E, \Gamma)$ forms a measurable field of Hilbert spaces where \[\Gamma:=\{s: P \to E: \textit{$s$ is measurable and}~ s(a) \in E(a), \forall a \in P\}.\]
\end{enumerate}

Set $\clh:=L^{2}(\Omega, \mathcal{F}, \mathbb{P})$. For $a \in P$, let $Q(a)$ be the projection on $\clh$ that corresponds to  the subspace $L^{2}(\Omega, \mathcal{F}_a, \mathbb{P})$. Then for $a \in P$ and $f \in \clh$, 
$Q(a)f=\mathbb{E}\{f|\mathcal{F}_a\}$. Note that 
\[
E:=\{(a,f) \in P \times L^{2}(\Omega, \mathcal{F}, \mathbb{P}): Q(a)f=f\}.\]
Thus to prove $(1)$ and $(3)$, it suffices to  show that $\{Q(a)\}_{a \in P}$ is a weakly measurable family of  projections, i.e. for $f,g \in \clh$, the map $P \ni a \to \langle Q(a)f|g \rangle$ is Borel measurable. 

\begin{lmma}
The family $\{Q(a)\}_{a \in P}$ is weakly measurable. 
\end{lmma}
\textit{Proof.} It suffices to prove that for $f,g$ in a total set, the map $P \ni a \to \langle Q(a)f|g \rangle \in \mathbb{C}$ is measurable. In view of Prop.\ref{total}, it is enough to show that for simple $L^{2}$-functions $u,v$ on $Y$, 
the map $P \ni a \to \langle Q(a)e^{\eta(u)}|e^{\eta(v)} \rangle \in \mathbb{C}$ is measurable. For $a \in P$, let $L_{a}:=X\backslash (X+a)$.

Write $u:=\sum_{i=1}^{n}a_{i}1_{B_i}$ and $v:=\sum_{i=1}^{n}b_{i}1_{B_i}$ with $\displaystyle \coprod_{i=1}^{n}B_i=Y$ and $\lambda(B_i)<\infty$. Calculate as follows to observe that 
\begin{align*}
Q(a)e^{\eta(u)}&=\mathbb{E}\{\prod_{i=1}e^{a_i \eta(B_i)}|\mathcal{F}_{a}\} \\
&=\mathbb{E}\{\prod_{i=1}^{n}e^{a_i \eta(B_i \cap L_a)}e^{a_i \eta(B_{i} \cap L_{a}^{c})}| \mathcal{F}_{a}\} \\
&=\prod_{i=1}^{n}e^{a_i \eta(B_i \cap L_a)}\prod_{i=1}^{n}\mathbb{E}\big(e^{a_i\eta(B_i \cap L_{a}^{c})}\big) ~~(\textrm{by the complete independence of $\eta$}).
\end{align*}
Similarly, $\displaystyle Q(a)e^{\eta(v)}=\prod_{i=1}^{n}e^{b_i\eta(B_i\cap L_a)}\prod_{i=1}^{n}\mathbb{E}\big(e^{b_i\eta(B_i \cap L_{a}^{c})}\big)$. Using the complete independence of the Poisson process , we arrive at the following equation
\[
\langle Q(a)e^{\eta(u)}|Q(a)e^{\eta(v)}\rangle=\prod_{i=1}^{n}\mathbb{E}\big(e^{(a_i+\overline{b_i})\eta(B_i \cap L_a)}\big)\prod_{i=1}^{n}\mathbb{E}\big(e^{a_i\eta(B_i \cap L_{a}^{c})}\big)\prod_{i=1}^{n}\mathbb{E}\big(e^{\overline{b_i}\eta(B_i \cap L_{a}^{c})}\big)\]

It is now sufficient to prove that for a measurable set $B$ of finite measure and a complex number $z$, the maps $P \ni a \to \mathbb{E}(e^{z\eta(B \cap L_{a})}) \in \mathbb{C}$ and $P \ni a \to \mathbb{E}(e^{z\eta(B \cap L_{a}^{c})})$ are measurable. Note that 
\[
\mathbb{E}(e^{z\eta(B \cap L_{a})})=exp(\lambda(B \cap L_a)(e^{z}-1))\] and
\[
\mathbb{E}(e^{z \eta(B \cap L_{a}^{c})})=exp(\lambda(B \cap L_{a}^{c})(e^{z}-1)).\]   By Tonelli's theorem, the  map  \[P \ni a \to \int 1_{B}(y)1_{X+a}(y)\lambda(dy) \in \mathbb{C}\] is measurable. Consequently,  for a measurable subset $B$ of finite measure, the maps $P \ni a \to \lambda(B \cap L_{a})\in \mathbb{C}$ and $P \ni a \to \lambda(B \cap L_{a}^{c})\in \mathbb{C}$ are measurable. 
The conclusion is now immediate and the proof is complete.  \hfill $\Box$

Let us fix a few notation. For $a \in P$, let $L_{-\infty,a}=Y \backslash (X+a)$ and  let $L_{a,\infty}=X+a$. For $a,b \in P$ with $a \leq b$, let $L_{a,b}=(X+a) \backslash (X+b)$. Note that for $a \leq b \leq c$, 
\begin{align*}
L_{a,c}&=L_{a,b} \sqcup L_{b,c} \\
Y&=L_{-\infty,a} \sqcup L_{a,b} \sqcup L_{b, \infty}.
\end{align*} 
For $a \in P$, we write $L_{a}=L_{0,a}$. 

\begin{lmma}
\label{measurability}
Let $B_1,B_2,\cdots, B_n$ be measurable subsets of $Y$ of finite measure and let $z_1,z_2,\cdots,z_n \in \mathbb{C}$ be given. Assume that $B_1, B_2, \cdots, B_n$ are disjoint. The maps 
 \begin{align*}
 ~~~P \ni &a \to e^{\sum_{i=1}^{n}z_i \eta(B_i \cap L_a)} \in L^{2}(\Omega, \mathcal{F},\mathbb{P})\\
 P \times P \ni &(a,b) \to e^{\sum_{i=1}^{n}z_i \eta((B_i-a) \cap L_b)} \in L^{2}(\Omega, \mathcal{F}, \mathbb{P}) \\
 P \times P &\ni (a,b) \to e^{\sum_{i=1}^{n}z_{i}\eta(B_i \cap L_{a+b,\infty})} \in L^{2}(\Omega, \mathcal{F}, \mathbb{P})
  \end{align*}
 are measurable. 
\end{lmma}
\textit{Proof.} Let $B$ be the complement of $\coprod_{i=1}^{n}B_i$ in $Y$. To prove that the first map is measurable, It suffices to prove that for a simple $L^{2}$-function $u$, the map \[P \ni a \to \mathbb{E}\Big(e^{\sum_{i=1}^{n}z_i \eta(B_i \cap L_a)}e^{\eta(u)}\Big) \in \mathbb{C}\] is measurable.
Let $u$ be a simple $L^{2}$-function. Write  $u=\sum_{j=1}^{m}w_{j}1_{C_j}$ with $\coprod_{j=1}^{m}C_{j}=Y$ and $\lambda(C_j)<\infty$. Set $A_{ij}=B_{i} \cap C_j$. 

A routine calculation using the complete independence of the Poisson process $\eta$ reveals that 
$
\mathbb{E}\Big(e^{\sum_{i=1}^{n}z_i\eta(B_i \cap L_a)}e^{\sum_{j=1}^{m}w_j\eta(C_j)}\Big)$ is the product of $\prod_{j=1}^{m}exp(\lambda(C_j \cap B \cap L_a)(e^{w_j}-1))$ and  \[\prod_{i,j}exp(\lambda(A_{ij} \cap L_a)(e^{z_i+w_j}-1))\prod_{j=1}^{m}exp(\lambda(C_j \cap L_{-\infty,0})(e^{w_j}-1))\prod_{j=1}^{m}exp(\lambda(C_j \cap L_{a,\infty})(e^{w_j}-1)).\] The measurability of each term follows from the fact that for a measurable set $B$ of finite measure the maps $P \ni a \to \lambda(B \cap L_a)$ and $P \ni a \to \lambda(B \cap L_{a,\infty})$ are measurable. The proofs of other assertions are similar and we omit the proof. \hfill $\Box$

\begin{lmma}
The map $E \times E \ni ((a,f),(b,g)) \to (a+b,fS_ag) \in E$ is measurable. 
\end{lmma}
\textit{Proof.} It suffices to show that for a simple function $L^{2}$-function $u$, \[((a,f),(b,g))\to \mathbb{E}(fS_age^{\eta(u)})\] is measurable. Write $u=\sum_{i=1}^{n}z_i1_{B_i}$ with $\coprod_{i=1}^{n}B_i=Y$ and $\lambda(B_i)<\infty$. 
Using the complete independence and the stationarity of the Poisson process $\eta$, note that $\mathbb{E}(fS_age^{\eta(u)})$ is 
\[
\mathbb{E}\Big(e^{\sum_i z_i \eta(B_i \cap L_{-\infty,0})}\Big)\mathbb{E}\Big(fe^{\sum_i z_i \eta(B_i \cap L_{a})}\Big)\mathbb{E}\Big(ge^{\sum_i z_i \eta((B_i-a) \cap L_b)}\Big)\mathbb{E}\Big(e^{\sum_i z_i \eta(B_i \cap L_{a+b,\infty})}\Big).\]
The required measurability conclusion follows from Lemma \ref{measurability}. This completes the proof. \hfill $\Box$

In short, we have the following theorem.

\begin{thm}
The set $E$, defined by Eq. \ref{product system}, together with its measurable and product structure is a product system over $P$. 
\end{thm}

\section{Identification of the product system associated to a stationary compound Poisson process}
In this section, we define a  product system associated to a stationary  compound Poisson process and identify it explicitly. First we consider the Poisson case.  We use the notation introduced in Section 3. For $a \in P$, let $\mathbbm{1}_{a} \in L^{2}(\Omega, \mathcal{F}_{a},\mathbb{P})$ be the constant function $1$.
Note that $(\mathbbm{1}_{a})_{a \in P}$ is a unit of $E$ (i.e. it is  a non-zero multiplicative measurable cross section of $E$).  For the rest of this paper, we fix a measurable logarithm, i.e. a measurable map $\ell:\mathbb{C}\backslash \{0\} \to \mathbb{C}$ such that 
$e^{\ell(z)}=z$ for every $z \in \mathbb{C} \backslash \{0\}$ and $\ell(x)=\log(x)$ if $x>0$. 

Let $\mathcal{S}$ denote the set of all complex valued simple  functions $f$ on $X$ which are square integrable (which is the same as integrable) such that $f(x) \neq -1$ for every $x \in X$. For $f \in \mathcal{S}$, set 
\begin{equation}
\label{u from f}
u_{f}(x):=\ell(1+f(x)).
\end{equation}
Note that $u_{f}$ is simple and square integrable. For $f \in \mathcal{S}$, let $\Sigma_{f}$ be the random variable defined by
\[
\Sigma_{f}:=\frac{e^{\eta(u_f)}}{\mathbb{E}\big(e^{\eta(u_f)}\big)}.\]

Observe that $\mathcal{S}$ is a dense subset of $L^{2}(X,\lambda)$. A simple calculation using Lemma \ref{Master equation} reveals that for $f,g \in \mathcal{S}$, 
\begin{equation}
\label{crucial equality 0}
\mathbb{E}(\Sigma_f \overline{\Sigma_g})=exp\Big(\int f(x)\overline{g(x)}\lambda(dx)\Big).
\end{equation}
Hence for $f,g \in \mathcal{S}$, 
\begin{equation}
\label{crucial equality}
\mathbb{E}\Big((\Sigma_f-\Sigma_g)\overline{(\Sigma_f-\Sigma_g)}\Big)=exp(\langle f|f \rangle)+exp(\langle g|g \rangle)-2 Re~ exp(\langle f|g \rangle).
\end{equation}
Let $f \in L^{2}(X,\lambda)$ be given. Choose a sequence $(f_n) \in \mathcal{S}$ such that $f_n \to f$ in $L^{2}(X,\lambda)$. Equation \ref{crucial equality} implies that $\Sigma_{f_n}$ converges in $L^{2}(\Omega,\mathcal{F},\mathbb{P})$. 
Making use of Equation \ref{crucial equality} again, it follows that the limit $\displaystyle \lim_{n \to \infty}\Sigma_{f_n}$ is independent of the chosen sequence $(f_n)$. Define 
\[
\Sigma_f=\lim_{n \to \infty} \Sigma_{f_n}.\]
Note that Equation \ref{crucial equality 0} is valid for $f,g \in L^{2}(X,\lambda)$. Observe  that for $a \in P$, if $f \in L^{2}(X\backslash (X+a))$ then $\Sigma_{f} \in L^{2}(\Omega,\mathcal{F}_{a}, \mathbb{P})$.

Consider the Hilbert space $H:=L^{2}(X,\lambda)$. For $a \in P$, let $V_{a}$ be the isometry on $H$ defined by the equation
\begin{equation}
\label{isometries}
V_{a}(f)(y):=\begin{cases}
 f(y-a)  & \mbox{ if
} y -a \in X,\cr
   &\cr
    0 &  \mbox{ if } y-a \notin X.
         \end{cases}
\end{equation}
Then $(V_a)_{a \in P}$ is a weakly measurable semigroup of isometries and hence is strongly continuous. %For $\phi \in L^{\infty}(X)$, we denote the multiplication operator on $H$ corresponding to $\phi$ by $M(\phi)$. 

\begin{ppsn}
\label{Poisson case}
The product system $E$ is isomorphic to the product system of the CCR flow associated to the isometric representation $V$. The map 
\[
L^{2}(X\backslash(X+a)) \ni e(f) \to \Sigma_{f} \in E(a)\]
for $a \in P$ provides an isomorphism between the product system associated to $V$ and $E$.  Here $e(.)$ denotes the exponential vectors.
\end{ppsn}
\textit{Proof.} Let $a \in P$ be given. For $f,g \in L^{2}(X\backslash(X+a))$, \[\mathbb{E}(\Sigma_f \overline{\Sigma_g})=exp(\langle f|g \rangle)=\langle e(f)|e(g)\rangle.\] Moreover the set $\{\Sigma_f: f \in L^{2}(X\backslash(X+a))\}$ is total in $E(a)$. This is because if $u$ is a non-negative simple function and if we set $f(x)=e^{u(x)}-1$ then $u_f=u$. Also the set $\{e(f): f \in L^{2}(X\backslash(X+a))\}$ is total in $\Gamma(L^{2}(X\backslash(X+a)))$. Consequently, there exists a unitary $\theta_a: \Gamma(L^{2}(X\backslash(X+a))) \to E(a)$ such that $\theta_a(e(f))=\Sigma_f$. Set $\displaystyle \theta:=\coprod_{a \in P}\theta_a$.

Let $a,b \in P$ and let $f \in L^{2}(X\backslash(X+a))$, $g \in L^{2}(X\backslash(X+b))$ be given. Observe that 
\[
(a,\Sigma_f)(b,\Sigma_g)=(a+b,\Sigma_h)
\]
where $h=f+V_a g$. Thus $\theta$ preserves the product structure. We leave the verification that $\theta$ is measurable to the reader. This completes the proof. \hfill $\Box$

The following is an immediate corollary of Prop. \ref{Poisson case} and Prop. 2.1 of \cite{Injectivity}.
\begin{crlre}
\label{decomposable vectors}
For $a \in P$, let $D(a)$ denote the set of decomposable vectors of $E(a)$. Then 
\[
D(a)=\{\lambda \Sigma_f: \lambda \in \mathbb{C}\backslash\{0\}, f \in L^{2}(X\backslash(X+a))\}.\]
\end{crlre}

\begin{rmrk}
Prop. \ref{Poisson case} could alternatively be derived by first proving Corollary \ref{decomposable vectors} by  imitating the proof of Prop. 2.1 of \cite{Injectivity}. Then it is clear that $E$ is decomposable and admits $(\mathbbm{1}_{a})$ as a unit. Appealing to Theorem 4.4 of \cite{Injectivity} yields Prop. \ref{Poisson case}. This  provides a conceptual explanation for the fact that product systems 
associated to stationary Poisson processes are CCR flows. 
\end{rmrk}

%\section{Identification of the product system associated to a compound Poisson process}
Next we associate a product system to a stationary compound Poisson process. First we collect a few preliminaries on compound Poisson processes from \cite{Last}. 
Let $(\mathbb{X},\mathcal{B})$ be a measurable space. We assume that $\mathcal{B}$ is countably generated. Denote the set of all $\sigma$-finite measures on $\mathbb{X}$ by $M(\mathbb{X})$. We endow $M(\mathbb{X})$ with the smallest $\sigma$-algebra which makes the map $\mu \to \mu(B)$ measurable for every $B \in \mathcal{B}$. A measurable mapping $\omega \to \xi(\omega) \in M(\mathbb{X})$ is called a random measure on $\mathbb{X}$ where $(\Omega,\mathcal{F},\mathbb{P})$ is an underlying probability space. Just like in the case of point processes, for every $\omega$, $\xi(\omega)$ is a measure and for every $B \in \mathcal{B}$, $\xi(B)$ is a random variable. 

The random measure that we will be  interested in are compound Poisson processes. Let $(Y,\mathcal{B}_{Y})$ be a measurable space and assume that $B_{Y}$ is countably generated. Suppose $\rho_0$ is a $\sigma$-finite measure on $Y$. Let $\nu$ be a ``L\'{e}vy measure" on $(0,\infty)$, i.e. the integral  $\int (r \wedge 1) \nu (dr)<\infty$ which is also equivalent to the fact that $\int (1-e^{-tr})\nu(dr)<\infty$ for every $t>0$. Let $\eta$ be a Poisson process on $Y \times (0,\infty)$ with intensity measure $\lambda:=\rho_0 \otimes \nu$. 

For $B \in B_{Y}$, let 
\[
\xi(B)=\eta(r1_{B}(y))=\int r1_{B}(y)\eta(d(y,r)).\]
Then $\xi$ is a random measure on $Y$. The random measure $\xi$ is called the $\rho_0$-symmetric compound Poisson process on $Y$ with L\'{e}vy measure $\nu$. The following are the basic properties of the compound Poisson process $\xi$. 
\begin{enumerate}
\item[(1)] The compound Poisson process $\xi$ is completely independent, i.e. for disjoint measurable subsets $B_1,B_2,\cdots,B_n$, the random variables $\xi(B_1),\xi(B_2),\cdots,\xi(B_n)$ are independent. 
\item[(2)] For $B \in B_{Y}$, the Laplace transform of the random variable $\xi(B)$ is given by
\[
\mathbb{E}(e^{-t\xi(B)})=exp\Big(-\rho_0(B)\int (1-e^{-tr})\nu(dr)\Big).\]
\end{enumerate}

Just like in the Poisson case, we can associate a product system to a stationary compound Poisson process. The setting is as before. Let $(Y,\mathcal{B}_{Y})$ be a measurable space on which $\mathbb{R}^{d}$ acts and let $\rho_0$ be a $\sigma$-finite measure on $Y$ which is $\mathbb{R}^{d}$ invariant. Suppose $X \subset Y$ is a measurable subset which is $P$-invariant. We also assume that the action of $P$ on $X$ is pure.  Let $\nu$ be a L\'{e}vy measure on $(0,\infty)$ and $\eta$ be the Poisson process on $Y \times (0,\infty)$ with intensity measure $\lambda:=\rho_0 \otimes \nu$.  Let $\xi$ be the $\rho_0$-symmetric compound Poisson process on $Y$ with L\'{e}vy measure $\nu$. Since $\rho_0$ is $\mathbb{R}^{d}$-invariant, the random measure $\xi$ is stationary, i.e. for a measurable subset $B$ of $Y$, $\xi(B+x)$ and $\xi(B)$ have the same distribution for every $x \in \mathbb{R}^{d}$. 

The action of $\mathbb{R}^{d}$ on $Y$ induces an action of $\mathbb{R}^{d}$ on $\widetilde{Y}:=Y \times (0,\infty)$ where the action is on the first coordinate. Let $X$ be a $P$-invariant measurable subset of $Y$ and set $\widetilde{X}:= X \times (0,\infty)$. Consider a probability space $(\Omega, \mathcal{F},\mathbb{P})$ which realises the Poisson process $\eta$. The $\sigma$-algebras $\mathcal{F}_{a}$, $\mathcal{F}_{a,b}$ etc.. are defined as in the Poisson case with $\eta$ replaced by $\xi$. 
The $\sigma$-algebras corresponding to the Poisson process $\eta$ associated to the data $(\widetilde{Y},\widetilde{X},\lambda)$ are denoted by $\widetilde{\mathcal{F}_a}$, $\widetilde{\mathcal{F}_{a,b}}$ etc...

Set \[
E:=\{(a,f): a \in P, f \in L^{2}(\Omega,\mathcal{F}_{a}, \mathbb{P})\}.\] The product rule on $E$ is defined exactly as in the Poisson case. Using the complete independence of $\xi$ and the stationarity of $\xi$, it is quite routine to check that the algebraic requirements for $E$ to be a product system is satisfied. It is possible to prove as  in the Poisson case that $E$ satisfies the measurability requirements. Howeover, it is automatic that the measurability requirements are met given that the Poisson case is already verified. We explain this below.

Let $(E,p)$ be a product system over $P$. Let $F \subset E$ be given. For $x \in P$, let $F(x)=E(x) \cap F$. We say that $F$ is a subsystem of $E$ if 
\begin{enumerate}
\item[(1)] for every $x \in P$, $F(x)$ is a non-zero closed subspace of $E(x)$, and
\item[(2)] for $x,y \in P$, $F(x)F(y) \subset F(x+y)$ and is total in $F(x+y)$. 
\end{enumerate}
Suppose $E$ is a product system and $F \subset E$ is a subsystem. We prove that $F$ is a measurable subset of $E$ and with the measurable structure induced from $E$, $F$ is a product system on its own right.
Realise $E$ as a product system of an $E_0$-semigroup, say $\alpha:=\{\alpha_x\}_{x \in P}$ on $B(\clh)$. 

For $x \in P$, let $\theta_x: E(x) \to E(x)$ be the orthogonal projection onto $F(x)$. Let $p_x \in \alpha_x(B(\clh))^{'}$ be such that 
$\theta_x(T)=p_xT$. Then $p_x$ is a non-zero projection for every $x \in P$. Note that for $u \in E(x)$ and $v \in E(y)$,
$\theta_{x+y}(uv)=\theta_x(u)\theta_y(v)$. This translates to the fact that $p_x \alpha_x(p_y)=p_{x+y}$ for $x,y \in P$. Then 
\[
F=\{(x,T): x \in P, T \in E(x), p_xT=T\}.\]
The required conclusion is immediate provided $(p_x)_{x \in P}$ is a weakly continuous family of projections. The latter assertion is proved  as in Prop. 8.9.9 of \cite{Arveson} with the aid of Theorem 10.8.1 of \cite{Hille}.

Let \[
\widetilde{E}:=\{(a,f): a \in P, f \in L^{2}(\Omega,\widetilde{\mathcal{F}}_{a},\mathbb{P})\}\] be the product system associated to the Poisson process $\eta$. Clearly $E$ is a subsystem of $\widetilde{E}$. Hence $E$ is a product system on its own right. 
Let $\widetilde{V}:=\{\widetilde{V}_{a}\}_{a \in P}$ be the semigroup of isometries on $L^{2}(\widetilde{X})$ induced by the action of $P$ on $\widetilde{X}$ (Eq. \ref{isometries}). 

Let us recall the Laplace functional of a Poisson process (Thm. 3.9, \cite{Last}). Suppose $u$ is a non-negative measurable function on $\widetilde{Y}$, then
\begin{equation}
\label{Laplace functional}
\mathbb{E}(e^{-\eta(u)})=exp\Big(-\int (1-e^{-u(y,r)})\lambda(d(y,r)\Big).
\end{equation}

\begin{thm}
The product systems $E$ and $\widetilde{E}$ coincide, i.e. $E=\widetilde{E}$. Consequently, $\widetilde{E}$ is isomorphic to a CCR flow. 
\end{thm}
\textit{Proof.} For $a \in P$, let $\theta_a: \widetilde{E}(a) \to \widetilde{E}(a)$ be the orthogonal projection corresponding to the subspace $E(a)$. For $a \in P$, $\theta_a(f)=\mathbb{E}\{f|\mathcal{F}_a\}$. Hence $\theta_{a}(\mathbbm{1}_a)=\mathbbm{1}_{a}$. 
Since $E$ is a subsystem of $\widetilde{E}$, it follows that 
$\{\theta_a\}_{a \in P}$ is a local projective cocycle of $\widetilde{E}$. But $\widetilde{E}$ is a CCR flow.
%The local projective cocycles of a CCR flow were computed in \cite{Injectivity}. 
Thanks to Prop. 5.1.1 of \cite{Injectivity}, Prop. \ref{Poisson case} and the fact that $\theta_a(\mathbbm{1}_a)=\mathbbm{1}_a$, it follows that there exists a projection $Q$ in the commutant of $\{\widetilde{V}_a: a \in P\}$ such that for $a \in P$, $f \in L^{2}(\widetilde{X}\backslash(\widetilde{X}+a))$, 
\[
\theta_{a}(\Sigma_f)=\Sigma_{Qf}.\]
The proof will be complete provided we show $Q=1$. Fix $a \in P$. It suffices to show that \[\mathbb{F}:=\{g \in L^{2}(\widetilde{X}\backslash(\widetilde{X}+a)): \Sigma_{g} \in L^{2}(\Omega,\mathcal{F}_{a},\mathbb{P})\}\] is total in $L^{2}(\widetilde{X}\backslash(\widetilde{X}+a))$ (For, then $H_a:=L^{2}(\widetilde{X}\backslash(\widetilde{X}+a)) \subset Ran(Q)$ for every $a \in P$ and the family of Hilbert spaces $(H_a)$ exhaust $L^{2}(\widetilde{X})$). 

Let $\mathcal{T}$ denote the set of all real valued measurable functions $f$ on $\widetilde{X}\backslash(\widetilde{X}+a)$ such that $-1<f(x,r) \leq 0$. For $f \in \mathcal{T}$, set $u_{f}(x,r)=\ell(1+f(x,r))=\log(1+f(x,r))$ and $v_f = -u_f$. Note that $v_f \geq  0$ for $f \in \mathcal{T}$. Let $B \subset  X \backslash(X+a)$ be a measurable subset of finite $\rho_0$-measure and let $c>0$ be given. Set $u:=c1_{B}$ and $\widetilde{u}(y,r)=ru(y)$. Note that $\frac{e^{-\xi(u)}}{\mathbb{E}(e^{-\xi(u)})}$ is a decomposable vector of $E(a)$ (also a decomposable vector of $\widetilde{E}(a))$ and has expectation $1$. Thus, by Corollary. \ref{decomposable vectors}, it follows that  there exists $g \in L^{2}(\widetilde{X}\backslash(\widetilde{X}+a))$ such that $\Sigma_{g}=\frac{e^{-\xi(u)}}{\mathbb{E}(e^{-\xi(u)})}$. 

For $f \in \mathcal{T}$, $\mathbb{E}(\Sigma_g\Sigma_f)=exp(\langle g|f \rangle)$. However 
\[
\mathbb{E}(\Sigma_g\Sigma_f)=\frac{\mathbb{E}(e^{-\eta(\widetilde{u}+v_f)})}{\mathbb{E}(e^{-\eta(\widetilde{u})})\mathbb{E}(e^{-\eta(v_f)})}.\]
 Applying Eq. \ref{Laplace functional} and simplifying, we obtain
\[
exp(\langle g|f\rangle)=\mathbb{E}(\Sigma_g\Sigma_f)=exp\Big(\int (-(1-e^{-ru(y)}))f(y,r) \lambda(d(y,r))\Big)=exp(\langle g_0|f \rangle)\]
where $g_0(y,r)=-(1-e^{-ru(y)})=-(1-e^{-cr})1_{B}(y)$.
Replacing $f$ by $tf$ for $t \in (0,1)$, we get $exp(t \langle g|f \rangle)=exp(t \langle g_0| f \rangle)$ for $t \in (0,1)$. Hence $\langle g|f \rangle=\langle g_0 |f \rangle$. But $\mathcal{T}$ is total in $L^{2}(\widetilde{X}\backslash(\widetilde{X}+a))$. Hence 
$g(y,r)=-(1-e^{-cr})1_{B}(y)$. 

Thus $\mathbb{F}$ contains the family of functions $\{-(1-e^{-cr})1_{B}(y): \rho_0(B)<\infty, c >0\}$. The totality of $\mathbb{F}$ in $L^{2}(\widetilde{X}\backslash(\widetilde{X}+a))$ follows from the fact that $\{1-e^{-cr}: c>0\}$ is total in $L^{2}((0,\infty),\nu)$. Hence $Q=1$ and the proof is now complete. \hfill $\Box$

\nocite{Tsirelson}
\nocite{Tsirelson_gauge}
\nocite{Tsi}
\nocite{Liebscher}
%\nocite{Arveson}
%\nocite{Kingman}
%\nocite{Arv}
%\nocite{Hilgert_Neeb}
%\nocite{Injectivity}
%\nocite{Murugan_Sundar_continuous}
\bibliography{references}

\def\cprime{$'$} \def\cprime{$'$} \def\cprime{$'$}
\providecommand{\bysame}{\leavevmode\hbox to3em{\hrulefill}\thinspace}
\providecommand{\MR}{\relax\ifhmode\unskip\space\fi MR }
\providecommand{\MRhref}[2]{%
  \href{http://www.ams.org/mathscinet-getitem?mr=#1}{#2}
}
\providecommand{\href}[2]{#2}
\begin{thebibliography}{10}

\bibitem{Arv_Fock4}
William Arveson, \emph{Continuous analogues of {F}ock space. {IV}. {E}ssential
  states}, Acta Math. \textbf{164} (1990), no.~3-4, 265--300.

\bibitem{Arveson}
\bysame, \emph{Noncommutative dynamics and {$E$}-semigroups}, Springer
  Monographs in Mathematics, Springer-Verlag, New York, 2003.

\bibitem{Arv}
\bysame, \emph{On the existence of {$E_0$}-semigroups}, Infin. Dimens. Anal.
  Quantum Probab. Relat. Top. \textbf{9} (2006), no.~2, 315--320.

\bibitem{Hille}
Einar Hille and Ralph~S. Phillips, \emph{Functional analysis and semi-groups},
  American Mathematical Society, Providence, R. I., 1974, Third printing of the
  revised edition of 1957, American Mathematical Society Colloquium
  Publications, Vol. XXXI.

\bibitem{Kingman}
J.~F.~C. Kingman, \emph{Poisson processes}, Oxford Studies in Probability,
  vol.~3, The Clarendon Press, Oxford University Press, New York, 1993, Oxford
  Science Publications.

\bibitem{Last}
G\"{u}nter Last and Mathew Penrose, \emph{Lectures on the {P}oisson process},
  Institute of Mathematical Statistics Textbooks, vol.~7, Cambridge University
  Press, Cambridge, 2018.

\bibitem{Liebscher}
Volkmar Liebscher, \emph{Random sets and invariants for (type {II}) continuous
  tensor product systems of {H}ilbert spaces}, Mem. Amer. Math. Soc.
  \textbf{199} (2009), no.~930.

\bibitem{Mark}
Daniel Markiewicz, \emph{On the product system of a completely positive
  semigroup}, J. Funct. Anal. \textbf{200} (2003), no.~1, 237--280.

\bibitem{Murugan_Sundar_continuous}
S.~P. Murugan and S.~Sundar, \emph{On the existence of {$E_0$}-semigroups---the
  multiparameter case}, Infin. Dimens. Anal. Quantum Probab. Relat. Top.
  \textbf{21} (2018), no.~2, 1850007, 20.

\bibitem{Injectivity}
S.~Sundar, \emph{Arveson's characterisation of {CCR} flows: the multiparameter
  case}, arxiv/math.OA:1906:05493v2.

\bibitem{Sundar_Notes}
S.~Sundar, \emph{Notes on {$E_0$}-semigroups}, available online at
  www.imsc.res.in/{$\sim$}ssundar.

\bibitem{Tsi}
Boris. Tsirelson, \emph{From random sets to continuous tensor products: answers
  to three questions of {W.} {A}rveson}, arxiv/math.FA:0001070.

\bibitem{Tsirelson}
Boris Tsirelson, \emph{Non-isomorphic product systems}, Advances in quantum
  dynamics, Contemp. Math., vol. 335, 2003, pp.~273--328.

\bibitem{Tsirelson_gauge}
\bysame, \emph{On automorphisms of type {II} {A}rveson systems (probabilistic
  approach)}, New York J. Math. \textbf{14} (2008), 539--576.

\end{thebibliography}
 \bibliographystyle{amsplain}

\noindent
{\sc S. Sundar}
(\texttt{sundarsobers@gmail.com})\\
         {\footnotesize  Institute of Mathematical Sciences (HBNI), CIT Campus, \\
Taramani, Chennai, 600113, Tamilnadu, INDIA.}\\

\end{document}